\documentclass[a4paper,12pt]{amsart}
\usepackage{amsfonts}
\usepackage{amssymb}
\usepackage{ifthen}
\usepackage{graphicx}
\usepackage{color}
\nonstopmode \numberwithin{equation}{section}
\newtheorem{thm}{Theorem}
\newtheorem{lem}{Lemma}
\newtheorem{cor}{Corollary}
\newtheorem{prop}{Proposition}

\newtheorem{conj}{Conjecture}

\theoremstyle{definition}
\newtheorem{defn}{Definition}
\newtheorem{example}{Example}

\newtheorem{ques}{Question}
\newtheorem{rem}{Remark}

\newtheorem{rems}{Remarks}

\newcounter {own}
\def\theown {\thesection  .\arabic{own}}

\newenvironment{pf}[1][]{%
 \vskip 3mm
 \noindent
 \ifthenelse{\equal{#1}{}}%
  {{\slshape Proof. }}%
  {{\slshape #1.} }%
 }%
{\qed\bigskip}

\DeclareMathOperator*{\esssup}{ess\,sup}

\newcounter{alphabet}
\newcounter{tmp}

\newcommand{\R}{{\mathbb R}}
\newcommand{\N}{{\mathbb N}}
\newcommand{\C}{{\mathbb C}}
\newcommand{\D}{{\mathbb D}}

\newcommand{\T}{{\mathbb{T}}}

\def\be{\begin{equation}}
\def\ee{\end{equation}}

\newcommand{\bee}{\begin{enumerate}}
\newcommand{\eee}{\end{enumerate}}

\newcommand{\blem}{\begin{lem}}
\newcommand{\elem}{\end{lem}}
\newcommand{\bthm}{\begin{thm}}
\newcommand{\ethm}{\end{thm}}
\newcommand{\bcor}{\begin{cor}}
\newcommand{\ecor}{\end{cor}}
\newcommand{\beg}{\begin{example}}
\newcommand{\eeg}{\end{example}}
\newcommand{\begs}{\begin{examples}}
\newcommand{\eegs}{\end{examples}}
\newcommand{\bdefn}{\begin{defn}}
\newcommand{\edefn}{\end{defn}}
\newcommand{\bprob}{\begin{prob}}
\newcommand{\eprob}{\end{prob}}
\newcommand{\bei}{\begin{itemize}}
\newcommand{\eei}{\end{itemize}}
\newcommand{\bqn}{\begin{ques}}
\newcommand{\eqn}{\end{ques}}
\newcommand{\bcon}{\begin{conj}}
\newcommand{\econ}{\end{conj}}
\newcommand{\bcons}{\begin{conjs}}
\newcommand{\econs}{\end{conjs}}
\newcommand{\bprop}{\begin{prop}}
\newcommand{\eprop}{\end{prop}}
\newcommand{\brem}{\begin{rem}}
\newcommand{\erem}{\end{rem}}
\newcommand{\brems}{\begin{rems}}
\newcommand{\erems}{\end{rems}}
\newcommand{\bo}{\begin{obser}}
\newcommand{\eo}{\end{obser}}
\newcommand{\bos}{\begin{obsers}}
\newcommand{\eos}{\end{obsers}}
\newcommand{\bpf}{\begin{pf}}
\newcommand{\epf}{\end{pf}}
\newcommand{\ba}{\begin{array}}
\newcommand{\ea}{\end{array}}
\newcommand{\beq}{\begin{eqnarray}}
\newcommand{\beqq}{\begin{eqnarray*}}
\newcommand{\eeq}{\end{eqnarray}}
\newcommand{\eeqq}{\end{eqnarray*}}

\newcommand{\ra}{\rightarrow}

\newcounter{minutes}
\divide\time by 60
\newcounter{hours}
\multiply\time by 60 \addtocounter{minutes}{-\time}

\begin{document}
\bibliographystyle{amsplain}
\title[Weighted composition operators on the class of subordinate functions]
{Weighted composition operators on the class of subordinate functions}

%

\author{Perumal Muthukumar}
\address{P. Muthukumar, Stat-Math Unit,
Indian Statistical Institute (ISI), Chennai Centre,
110, Nelson Manickam Road,
Aminjikarai, Chennai, 600 029, India.}
\email{pmuthumaths@gmail.com}

\author{Saminathan Ponnusamy  
}
\address{S. Ponnusamy, Department of Mathematics,
Indian Institute of Technology Madras, Chennai-600 036, India.}
\email{samy@iitm.ac.in, samy@isichennai.res.in }

\


\begin{abstract}
In this article, we study the weighted composition operators preserving the class $\mathcal{P}_{\alpha}$ of analytic functions subordinate to
$\frac{1+\alpha z}{1-z}$ for $|\alpha|\leq 1, \alpha \neq -1$. Some of its consequences and examples for some special cases are presented.
\end{abstract}
\thanks{
File:~\jobname .tex,
          printed: \number\day-\number\month-\number\year,
          \thehours.\ifnum\theminutes<10{0}\fi\theminutes
}
\maketitle
\pagestyle{myheadings}
\markboth{P. Muthukumar and S. Ponnusamy}{Weighted composition operators on the class of subordinate functions}

 \section{Introduction}
Let $\mathcal{H}(\mathbb{D})$ denote the class of analytic functions defined on the unit disk $\D= \{z \in \C:\, |z|<1\}$
with the (metrizable) topology of uniform convergence on compact subsets of $\D$ and we denote the boundary of $\D$ by $\T$.
Weighted composition operator is a combination of multiplication and composition operators.
These operators are mainly studied in various Banach spaces or Hilbert spaces of $\mathcal{H}(\mathbb{D})$.

Recently,
Ar\'evalo et al. \cite{arxiv W-CO}
initiated the study of weighted composition operator restricted to the Carath\'eodory class $\mathcal{P}_{1}$, which consists of all
$f \in \mathcal{H}(\mathbb{D})$ with positive real part and with a normalization $f(0)=1$. Clearly the class $\mathcal{P}_{1}$
is  not a linear space but it is helpful to solve some extremal problems in geometric
function theory. See \cite{Hallenbeck:Book}.

In this article, we generalize the recent work of Ar\'evalo et al. \cite{arxiv W-CO} by considering
weighted composition operators preserving the class $\mathcal{P}_{\alpha}$. This class is connected
with various geometric subclasses of $\mathcal{H}(\mathbb{D})$ in the univalent function theory
(see \cite{Duren,Hallenbeck:Book,Pommerenke:Book}).
Since the class $\mathcal{P}_{\alpha}$ is not a
linear space, for a given map on  $\mathcal{P}_{\alpha}$, questions about operator theoretic properties are not
meaningful. However, one can talk about special classes of self-maps of $\mathcal{P}_{\alpha}$ and fixed points
of those maps. This is the main purpose of this article.

The article is organized as follows. In Section \ref{MP4Sec2}, we introduce the class $\mathcal{P}_{\alpha}$ and
list down some basic properties about
this class. In Section \ref{MP4Sec3}, we give characterization for weighted composition operators to be self-maps
of the class $\mathcal{P}_{\alpha}$ (see Theorem \ref{MP4Main1}). The above situation is analyzed closely for various
special cases of symbols in Section \ref{MP4Sec4}. In Section \ref{MP4Sec5}, we present some simple examples.

\section{Some preliminaries about the class $\mathcal{P}_{\alpha}$}\label{MP4Sec2}
 For $f$ and $g \in \mathcal{H}(\mathbb{D})$, we say that \textit{$f$ is subordinate to $g$} (denoted by $f(z)\prec g(z)$ or $f\prec g$) if
there exists an analytic function $\omega:\D\rightarrow \D$ such that $\omega(0)=0$ and $f=g\circ \omega$. If $f(z)\prec z$, then
$f$ is called \textit{Schwarz} function. For $|\alpha|\leq 1, \alpha \neq -1$, define $h_\alpha$ on $\D$ by
$ h_\alpha(z)=\frac{1+\alpha z}{1-z}$ and the half plane $\mathbb{H}_\alpha$  is described by
$$ \mathbb{H}_\alpha :=
h_\alpha(\D)=\{ w\in \C: 2{\rm Re\,}\{(1+\overline{\alpha})w\}>1-|\alpha|^2 \}.
$$
In particular, if $\alpha \in \R$ and $-1<\alpha\leq1$, then
$$
h_\alpha(\D)=\{ w\in \C: {\rm Re\,} w> (1-\alpha)/2\}
$$
so that ${\rm Re\,}h_\alpha(z)> (1-\alpha)/2$ in $\D$.

For $|\alpha|\leq 1, \alpha \neq -1$, it is natural to consider the class $\mathcal{P}_{\alpha}$ defined by
$$ \mathcal{P}_{\alpha}:= \{f \in \mathcal{H}(\mathbb{D}) : f(z)\prec  h_\alpha(z) \}.
$$
It is worth to note that for every $f \in \mathcal{P}_{\alpha}$, there is an unique Schwarz function $\omega$
such that
$$f(z)=\frac{1+\alpha \omega(z)}{1-\omega(z)}.$$
It is well-known \cite[Lemma 2.1]{Pommerenke:Book} that, if $g$ is an univalent analytic function on $\D$, then
$f(z)\prec g(z)$  if and only if $f(0)=g(0)$ and $f(\D)\subset g(\D)$. In view of this result, the class
$\mathcal{P}_{\alpha}$ can be stated in an equivalent form as
$$ \mathcal{P}_{\alpha}:= \{f \in \mathcal{H}(\mathbb{D}) :f(0)=1, f(\D)\subset \mathbb{H}_\alpha\}.
$$
We continue the discussion by stating a few basic and useful properties of the class $\mathcal{P}_{\alpha}$.

\bprop
Suppose $f \in  \mathcal{P}_{\alpha}$ and $f(z)=1+\sum_{n=1}^{\infty}a_n z^n$, then $|a_n|\leq |\alpha+1|$ for all $n\in \N$.
The bound is sharp as the function $h_\alpha(z)=1+\sum_{n=1}^{\infty}(1+\alpha)z^n$ shows.
\eprop
\bpf
This result is an immediate consequence of Rogosinski's result \cite[Theorem X]{subordination} (see also \cite[Theorem 6.4(i), p.~195]{Duren})
because $h_\alpha(z)$ (and hence, $(h_\alpha(z)-1)/(1+\alpha)$) is a convex function.
\epf

\bprop {\rm (Growth estimate)}
Let $f \in  \mathcal{P}_{\alpha}$. Then  for all $z\in \D$, one has
$$ \frac{1-|\alpha z|}{1+|z|}\leq |f(z)|\leq \frac{1+|\alpha z|}{1-|z|}.
$$
\eprop
\bpf
This result trivially follows from clever use of classical Schwarz lemma and the triangle inequality.
\epf

From the `growth estimate' and the familiar Montel's theorem on normal family, one can easily get the following result.

\bprop
The class $\mathcal{P}_{\alpha}$ is a compact family in the compact-open topology (that is, topology of uniform convergence
on compact subsets of $\D$).
\eprop

Because the half plane $\mathbb{H}_\alpha$ is convex, the following result is obvious.
\bprop
The class $\mathcal{P}_{\alpha}$ is a convex family.
\eprop


For $p\in(0,\infty)$, the Hardy space $H^{p}$ consists of analytic functions $f$ on $\mathbb{D}$ with
$$\|f\|_{p}:=\sup\limits_{ r \in [0,1)}\left \{\frac{1}{2\pi}\int_0^{2\pi} |f(re^{i\theta})|^{p}\,d\theta\right \}^{\frac{1}{p}}
$$
is finite and $H^\infty$ denotes the set of all bounded analytic functions on $\mathbb{D}$. We refer to
\cite{Duren:Hpspace} for the theory of Hardy spaces.
By Littlewood's subordination theorem \cite[Theorem 2]{littlewood thm}, it follows that if  $f\prec g$ and
$g\in H^p$ for some $0<p\leq \infty$, then $f\in H^p$ for the same $p$. As a consequence we easily have the following.

\bprop
The class $\mathcal{P}_{\alpha}$ is a subset of the Hardy space $H^p$ for each $0<p<1$.
\eprop
\bpf
Because $(1-z)^{-1} \in H^p$ for each $0<p<1$, it follows easily that $h_\alpha \in H^p$ for each $0<p<1$ and
for $|\alpha|\leq 1, \alpha \neq -1$. The desired conclusion follows.
\epf

\brem
Although $\mathcal{P}_{\alpha}$
does not posses the linear structure, due to being part of $H^p$,  the results on $H^p$ space, such as results about
boundary  behavior, are valid for functions in the class $\mathcal{P}_{\alpha}$
\erem

\section{Weighted composition on $\mathcal{P}_{\alpha}$}\label{MP4Sec3}

For an analytic self-map $\phi$ of $\D$, the composition operator $C_\phi$ is defined by
$$
C_\phi(f)=f\circ \phi \mbox{~for~}
f \in \mathcal{H}(\mathbb{D}). $$
One can refer \cite{Cowen:Book}, for the study of composition operators on various
function spaces on the  unit disk. Throughout this article, $\alpha$ denotes a complex number such that
$|\alpha|\leq 1, \alpha \neq -1$ , unless otherwise stated explicitly and $\phi$ denotes an analytic self-map of $\D$.
The following result deals with composition operator when it is restricted to the class
$\mathcal{P}_{\alpha}$.

\bprop
The composition operator $C_\phi$ induced by the symbol $\phi$, preserves the class $\mathcal{P}_{\alpha}$ if and
only if $\phi$ is a Schwarz function.
\eprop
\bpf
Suppose that $C_\phi$  preserves the class $\mathcal{P}_{\alpha}$. Then $C_\phi (h_\alpha) \in \mathcal{P}_{\alpha}$,
and thus
$$\frac{1+\alpha \phi(0)}{1-\phi(0)}=1.
$$
This gives that $\phi(0)=0$ which implies that $\phi$ is a Schwarz function.
The converse part holds trivially.
\epf

For a given analytic self-map $\phi$ of $\D$ and analytic map $\psi$ of $\D$, the corresponding weighted
composition operator $C_{\psi,\phi}$ is defined by
$$ C_{\psi,\phi}(f)= \psi (f\circ \phi) \mbox {~for~} f \in \mathcal{H}(\mathbb{D}).
$$

If $\psi\equiv 1$, then $C_{\psi,\phi}$ reduced to a composition operator $C_{\phi}$ and if $\phi(z)=z$ for all $z\in \D$, then $C_{\psi,\phi}$
reduced to a multiplication operator $M_{\psi}$.
For a given analytic map $\psi$ of $\D$, the corresponding multiplication operator $M_\psi$ is then defined by
$$ M_\psi(f)= \psi f \mbox {~for~} f \in \mathcal{H}(\mathbb{D}).
$$
The characterization of $M_{\psi}$ that preserves the class $\mathcal{P}_{\alpha}$ is given in Section \ref{MP4Sec4}.

Banach begun the study of weighted composition operators. In \cite{Banach}, Banach proved the classical Banach-Stone theorem, that is, the surjective isometries between the spaces of continuous real-valued functions on a closed and bounded interval are certain weighted composition operators. In \cite{Forelli}, Forelli proved that the isometric isomorphism of the Hardy space $H^p, (p\neq 2)$ are also weighted composition operators. The same result for the case of Bergman space is proved by Kolaski in \cite{Bergman isometry}.

The study of weighted composition operators can be viewed as a natural generalization of the well known field in the analytic function theory, namely, the composition operators. Moreover, weighted composition operators appear in applied areas such as dynamical systems and evolution equations. For example, classification of dichotomies in certain dynamical systems is connected to weighted composition operators, see \cite{Dynamical}.

In this section, we discuss
weighted composition operator that preserves $\mathcal{P}_{\alpha}$. Before, we do this, let us recall some useful results from the theory of extreme points.

\blem {\rm (\cite[Theorem 5.7]{Hallenbeck:Book})}
Extreme points of the class $\mathcal{P}_{\alpha}$  consists of functions given by
$$
f_\lambda(z)= \frac{1+\alpha \lambda z}{1-\lambda z},~~ |\lambda|=1.
$$
\elem

A point $p$ of a convex set $A$ is called \textit{extreme point} if $p$ is not a interior point of any line segment
which entirely lies in $A$. We denote, the set of all extreme points of the class $\mathcal{P}_{\alpha}$ by
$\mathcal{E}_{\alpha}$. That is, $\mathcal{E}_{\alpha}= \{f_\lambda: |\lambda|=1\}.$ Now, we recall a well-known
result by Krein and Milman \cite{Milman}.

\blem {\rm  (\cite[Theorem 4.4]{Hallenbeck:Book})}
Let $X$ be a locally convex, topological vector space and $A$ be a convex, compact subset of $X$. Then, the closed
convex hull of extreme points of $A$ is equal to $A$.
\elem

The original version of it is proved in \cite{Milman}.
On $\mathcal{H}(\mathbb{D})$, $f_n$ converges to $f$
denoted by $f_n \xrightarrow{\text{u.c}} f$. It is easy to see that $ C_{\psi,\phi}(f_n) \xrightarrow{\text{u.c}}
C_{\psi,\phi}(f)$ whenever $f_n \xrightarrow{\text{u.c}} f$. Thus, $ C_{\psi,\phi}$ is continuous on
$\mathcal{H}(\mathbb{D})$ (in particular on $\mathcal{P}_{\alpha}$).

\bprop
Suppose that $C_{\psi,\phi}$ preserves the class $\mathcal{P}_{\alpha}$. Then $\phi$ is a Schwarz function and
there exists a Schwarz function $\omega$ such that
$$\psi= h_\alpha \circ \omega= \frac{1+\alpha \omega}{1-\omega}.
$$
\eprop
\bpf
Suppose that $C_{\psi,\phi}$ preserves the class $\mathcal{P}_{\alpha}$. Take $f\equiv 1$ to be a constant function,
which belongs to $\mathcal{P}_{\alpha}$. Thus, $C_{\psi,\phi} (f)=\psi \in \mathcal{P}_{\alpha}$ and hence,
there exists a Schwarz function $\omega$ such that
$$\psi= h_\alpha \circ \omega= \frac{1+\alpha \omega}{1-\omega}.
$$
In particular, $\psi(0)=1$.

Since $h_\alpha \in \mathcal{P}_{\alpha}$, we have $\psi (0) (h_\alpha ( \phi(0)))=1$, which gives $\phi(0)=0$.
Hence $\phi$ will be a Schwarz function.
\epf

In view of above result, from now on, we will assume that $\psi= h_\alpha \circ \omega= \frac{1+\alpha \omega}{1-\omega}$
and $\phi,\,\omega$ are Schwarz functions.

\bthm\label{MP4Main1}
Let $\phi$, $\omega$ and $\psi$ be as above. Then, $C_{\psi,\phi}$ preserves the class $\mathcal{P}_{\alpha}$ if and only if
\begin{equation}\label{mp4 eq1}
 2Q(\omega) |\phi|< (1-|\omega|^2)+P(\omega) |\phi|^2 \mbox{~on~} \D,
\end{equation}
where, $P(\omega)= |\alpha \omega|^2-|1+(\alpha-1)\omega|^2$ and
$Q(\omega)= |(\alpha-1)|\omega|^2 + \overline{\omega}-\alpha \omega|$.
\ethm
\bpf
At first we prove that, $C_{\psi,\phi}$ preserves the class $\mathcal{P}_{\alpha}$ which is equivalent to
the inclusion $C_{\psi,\phi}(\mathcal{E}_{\alpha})\subset\mathcal{P}_{\alpha}$. To do this, we suppose that
$C_{\psi,\phi} (\mathcal{E}_{\alpha})\subset \mathcal{P}_{\alpha}$. Since $\mathcal{P}_{\alpha}$ is a
convex family, we obtain
$$ C_{\psi,\phi} (\mbox{~convex hull~}(\mathcal{E}_{\alpha}))\subset \mathcal{P}_{\alpha}.
$$
Now, by Krein-Milman theorem and the fact that $ C_{\psi,\phi}$ is continuous on  a compact family
$\mathcal{P}_{\alpha}$, we see that $C_{\psi,\phi} (\mathcal{P}_{\alpha})\subset \mathcal{P}_{\alpha}$.
The converse part is trivial.

Next, we prove that $C_{\psi,\phi} (\mathcal{E}_{\alpha})\subset \mathcal{P}_{\alpha}$ if and only if
$$ 2Q(\omega) |\phi|< (1-|\omega|^2)+P(\omega) |\phi|^2 \mbox{~on~} \D.
$$
Assume that $C_{\psi,\phi} (\mathcal{E}_{\alpha})\subset \mathcal{P}_{\alpha}$. This gives
$\psi (f_\lambda \circ \phi)\in \mathcal{P}_{\alpha}$ for all $|\lambda|=1$. Thus, for all $|\lambda|=1$,
there exists a Schwarz function $\omega_\lambda$ such that
$\psi (f_\lambda \circ \phi)= h_\alpha \circ \omega_\lambda$. That is,
$$ \frac{1+\alpha\omega}{1-\omega} \frac{1+\alpha\lambda\phi}{1-\lambda\phi}=
\frac{1+\alpha\omega_\lambda}{1-\omega_\lambda}.
$$
Solving this equation for $\omega_\lambda$, we get that
$$ \omega_\lambda=\frac{\omega+\lambda\phi+(\alpha-1)\lambda\omega\phi}{1+\alpha\lambda\phi\omega}.
$$

For each $|\lambda|=1$, $\omega_\lambda$ is a Schwarz function if and only if
$$ |\omega+\lambda\phi+(\alpha-1)\lambda\omega\phi|^2< |1+\alpha\lambda\phi\omega|^2
~\mbox{ for all }~ |\lambda|=1,
$$
which is equivalent to
$$ 2{\rm Re~}(\lambda\phi\{(\alpha-1)|\omega|^2 + \overline{\omega}-\alpha \omega\}) <
(1-|\omega|^2)+|\phi|^2(|\alpha \omega|^2-|1+(\alpha-1)\omega|^2),
$$
for all $|\lambda|=1$. By taking supremum over $\lambda$ on both sides, the last inequality
gives \eqref{mp4 eq1}.
The converse part follows by repeating the above arguments in the reverse direction.
 \epf

\brem \label{p and q}
Suppose that $\alpha=a+ib$ and $\omega(z)=u(z)+iv(z)$. Then,
\beqq
P(\omega)&=& |\alpha \omega|^2-|1+(\alpha-1)\omega|^2\\
&=& a(|\omega|^2-1)+(a-1)|\omega-1|^2+2bv.
\eeqq
Set $q(\omega)=(\alpha-1)|\omega|^2 + \overline{\omega}-\alpha \omega $ so that
$Q(\omega)= |q(\omega)|$. Upon simplifying, we get that
$$q(\omega)= (\alpha-1)\,\overline{\omega}\,(\omega-1)-2i\alpha\,v=(\alpha-1)( |\omega|^2- \omega) - 2iv$$ and thus
\begin{equation}\label{mp4 eq2}
 q(\omega)=[(a-1)(|\omega|^2-u)+bv]+i[b(|\omega|^2-u)-v(a+1)].
\end{equation}
Also, it is easy to see that
\begin{equation}\label{mp4 eq3}
-q(\omega)= |1-\omega|^2\psi+(|\omega|^2-1)  ~\mbox{ with }~ \psi=\frac{1+\alpha \omega}{1-\omega}.
\end{equation}
\erem

\section{Special cases}\label{MP4Sec4}

In this section, first we recall some familiar results on Hardy space $H^p$ which will help
the smooth traveling of this article.

\bprop {\rm (\cite[Theorem 1.3]{Duren:Hpspace})}\label{radial limit}
For every bounded analytic function $f$ on $\D$, the radial limit $\lim\limits_{r\rightarrow 1} f(re^{i\theta})$
exists almost everywhere (abbreviated by a.e.).
\eprop

In view of Proposition \ref{radial limit}, every Schwarz function has radial limit a.e. and using the fact that the function
$h_{\alpha}$ has radial limit a.e., it is easy to see that, every function $f\in \mathcal{P}_{\alpha}$ has radial limit
a.e.  on $\T$. Also, it is well-known that (see \cite[Section 2.3]{Duren:Hpspace})
$$ \sup\limits_{|z|<1}|f(z)|= \esssup_{0\leq \theta\leq 2 \pi}|f(e^{i\theta})|,
$$
for every $f\in H^\infty$. Now, we will state a classical theorem of Nevanlinna.

\bprop {\rm (\cite[Theorem 2.2]{Duren:Hpspace})}\label{posite mre-radial limit}
If $f\in H^p$ for some $p>0$ and its radial limit $f(e^{i\theta})=0$ on a set of positive measure, then $f\equiv 0$.
\eprop

Since every Schwarz function $f$ belongs to $H^\infty$ and every $f\in \mathcal{P}_{\alpha}$
belongs to $H^p$ for $0<p<1$, the above result is valid for functions in the class
$\mathcal{P}_{\alpha}$ and for Schwarz functions.

An analytic function $f$ on $\D$ is said to be an \textit{ inner function } if $|f(z)|\leq 1$ for all $z\in \D$
and its radial limit $|f(\zeta)|=1$ a.e. on $|\zeta|=1$.

\bthm\label{MP4Main2}
Suppose that $\phi$ and $\omega$ are Schwarz functions, $\phi$ is inner and $\psi=\frac{1+\alpha \omega}{1-\omega}$.
Then,  $C_{\psi,\phi}$ preserves the class $\mathcal{P}_{\alpha}$ if and only if
$\psi\equiv 1$ (i.e., $\omega\equiv 0$).
\ethm
\bpf
If $\psi\equiv 1$, then $C_{\psi,\phi}$ becomes a composition operator $C_{\phi}$ and thus,
$C_{\psi,\phi}$ preserves the class $\mathcal{P}_{\alpha}$, because $\phi$ is a Schwarz function.

Conversely, suppose that $C_{\psi,\phi}$ preserves the class $\mathcal{P}_{\alpha}$.
Then, by Theorem \ref{MP4Main1}, one has the inequality
$$
 2Q(\omega) |\phi|< (1-|\omega|^2)+P(\omega) |\phi|^2 \mbox{~on~} \D.
$$

With abuse of notation, we denote the radial limits of $\phi$, $\omega$ and $\psi$
again by $\phi$, $\omega$ and $\psi$, respectively. Also, let $\alpha=a+ib$ and $\omega(z)=u(z)+iv(z)$.
By allowing $|z|\ra 1$ in \eqref{mp4 eq1}, we get that
$$
2Q(\omega) \leq (1-|\omega|^2)+P(\omega)  \mbox{~a.e. on~} \T,
$$
which after computation is equivalent to
$$
Q(\omega) \leq (a-1)(|\omega|^2-u)+bv  \mbox{~a.e. on~} \T.
$$
In view of \eqref{mp4 eq2} in Remark \ref{p and q}, the above inequality can rewritten as
$$
|q(\omega)| \leq {\rm Re\,}[q(\omega)] \mbox{~a.e. on~} \T
$$
which gives that
${\rm Im\,}[q(\omega)]=0 \mbox{~a.e. on~} \T$.
Again by using \eqref{mp4 eq3} in Remark \ref{p and q},  we have
$$|1-\omega|^2{\rm Im}(\psi)=0 \mbox{~a.e. on~} \T.$$
Analyzing the function $\omega$ through the classical theorem of Nevanlinna (see Proposition \ref{posite mre-radial limit}),
one can get that ${\rm Im\,}\psi=0$ a.e. on $\T$. Now the proof of $\psi\equiv 1$ is as follows:

Consider the analytic map $f=e^{-i(\psi-1)}$. Then, $|f|=e^{{\rm Im\,}\psi}=1$ a.e. on $\T$ and
$$
1=f(0)\leq \sup\limits_{|z|<1}|f(z)|= \esssup_{0\leq \theta\leq 2 \pi}|f(e^{i\theta})|=1.
$$
Hence, by the maximum modulus principle, we get that $f\equiv 1$ which gives $\psi\equiv 1$.
\epf
%

\bcor
$M_\psi$ preserves the class $\mathcal{P}_{\alpha}$ if and only if
$\psi\equiv 1$.
\ecor
\bpf
The desired result follows if we set $\phi(z)\equiv z$ in Theorem \ref{MP4Main2}.
\epf

\bthm\label{MP4Main3}
Suppose that $\alpha$ is a real number, $\phi$ and $\omega$ are Schwarz functions, $\omega$ is an inner function
and $\psi=\frac{1+\alpha \omega}{1-\omega}$.
Then,  $C_{\psi,\phi}$ preserves the class $\mathcal{P}_{\alpha}$ if and only if
$\phi$ is identically zero.
\ethm
\bpf
If $\phi\equiv 0$, then $C_{\psi,\phi}$ becomes a constant map $\psi$ and
hence it preserves $\mathcal{P}_{\alpha}$.
Conversely, suppose that $C_{\psi,\phi}$ preserves the class $\mathcal{P}_{\alpha}$.
Then, by Theorem \ref{MP4Main1},
$$
 2Q(\omega) |\phi|< (1-|\omega|^2)+P(\omega) |\phi|^2 \mbox{~on~} \D.
$$
By allowing $|z|\ra 1$, we get that
$$
2|1-\omega|^2|\psi|\,|\phi|\leq (2 {\rm Im\,}\alpha\,{\rm Im\,}\omega+
({\rm Re\,}\alpha-1)|1-\omega|^2)|\phi|^2 \mbox{~a.e. on~} \T,
$$
from which we obtain that
$$
|1-\omega|^2|\psi|\,|\phi|\leq 0 \mbox{~a.e. on~} \T.
$$
By the hypothesis on $\omega$ and $\psi$, and the classical theorem of Nevanlinna, we find
that $\phi\equiv 0$.
\epf

Here is an easy consequence of Theorem \ref{MP4Main3}.

\bcor
Let $\alpha$ be a real number, $\phi$ and $\omega$ are Schwarz functions and that $\phi\not \equiv 0$.
Suppose that $C_{\psi,\phi}$ preserves the class $\mathcal{P}_{\alpha}$, and
$$ E=\{ \zeta \in \T: |\omega(\zeta)|=1\}.
$$
Then the Lebesgue arc length measure of the set $E$ is zero, i.e., $m(E)=0$.
\ecor

\section{Examples for special cases}\label{MP4Sec5}

In this section, we give specific examples of $\phi$ and $\psi$ so that
$C_{\psi,\phi}$ preserves the class $\mathcal{P}_{\alpha}$. For a bounded analytic function on $\D$,
we denote $\sup\limits_{|z|<1}|f(z)|$ by $\|f\|$.
\beg
Suppose that $\|\phi\|<1$. If $\|\omega\|<\frac{1-\|\phi\|}{1+\|\phi\|}$,
then $C_{\psi,\phi}$ preserves the class $\mathcal{P}_{\alpha}$, for $\alpha \in [0,1]$.
\eeg
\bpf In view of Theorem \ref{MP4Main1}, it suffices to verify the inequality \eqref{mp4 eq1}.
This inequality can be rewritten as
$$
 2|(1- \alpha)\,\overline{\omega}\,(\omega-1)+2i\alpha\, {\rm Im\,}\omega|\,|\phi|+(1- \alpha)(|\omega|^2 +|1-\omega|^2|\phi|^2-1)
 < \alpha(1-|\omega|^2)(1-|\phi|^2).
$$
We may set $\|\omega\|=A$ and $\|\phi\|=B$. Thus, it is enough to check that
$$
 2(1- \alpha)A(A+1)B+4\alpha AB+(1- \alpha)(A^2 +(1+A)^2B^2-1)
 < \alpha(1-A^2)(1-B^2)
$$
which is equivalent to
$$
[A+B+AB-1][(1- \alpha)(A+B+AB+1)+\alpha (A+B-AB+1)]<0.
$$
This yields the condition $A+B+AB-1<0$. This means that $A<\frac{1-B}{1+B}$ and the desired conclusion follows.
\epf

Since the condition $A+B+AB-1<0$ gives $B<\frac{1-A}{1+A}$, we have the following result.
\beg
Suppose that $\|\omega\|<1$. If $\|\phi\|<\frac{1-\|\omega\|}{1+\|\omega\|}$,
then $C_{\psi,\phi}$ preserves the class $\mathcal{P}_{\alpha}$, for $\alpha \in [0,1]$.
\eeg

\beg
Suppose that $\phi(z)=z(az+b)$, where $a$ and $b$ are non-zero real numbers such that $|a|+|b|=1$.
Take $\omega(z)=z(cz+d)$ with
$$c=-\frac{ab}{K} ~\mbox{ and } ~d=\frac{1-(a^2+b^2)}{K} ~\mbox{ for $K>2+\sqrt{5}$}.
$$
Then $C_{\psi,\phi}$ preserves the class $\mathcal{P}_{1}$.
\eeg
\bpf
Clearly $|\phi(z)|^2\leq a^2+ b^2+ 2abx$ for  $z=x+iy$ and thus,
$$
0< 1-(a^2+ b^2)-2abx\leq (1-|\phi|^2).
$$
Also note that
$$ |{\rm Im\,}\omega|\leq |2cx+d| = \frac{1-(a^2+ b^2)-2abx}{K}
$$
and
$$
|\omega(z)|\leq |c|+|d|= \frac{1-|ab|-(|a|-|b|)^2}{K}\leq \frac{1}{K}.
$$
The inequality \eqref{mp4 eq1} for $\alpha=1$ reduces to
$$
4|\phi|\,|{\rm Im\,}\omega|<(1-|\omega|^2)(1-|\phi|^2).
$$
Since $4|\phi|\,|{\rm Im\,}\omega|\leq 4|{\rm Im\,}\omega|\leq 4|2cx+d|$ and
$$ \left (1-\frac{1}{K^2}\right )K|2cx+d|\leq(1-|\omega|^2)(1-|\phi|^2),
$$ to verify the inequality \eqref{mp4 eq1},
it suffices to verify the inequality
$$ \frac{4}{K}< 1-\frac{1}{K^2}, ~\mbox{ i.e., }~ K^2-4K-1>0.
$$
This gives the condition $K>2+\sqrt{5}$ and the proof is complete.
\epf

\brem\label{Rem3-new}
By letting $\alpha=0$ in the Theorem \ref{MP4Main1}, we see that $C_{\psi,\phi}$ preserves the class
$\mathcal{P}_{0}$ if and only if $|1-\omega| \, |\phi|+|\omega|<1 $ on $\D.$
\erem
\beg
If $|\phi| \leq |\omega|<\sqrt{2}-1$ on $\D$, then $C_{\psi,\phi}$ preserves $\mathcal{P}_{0}$.
\eeg
\bpf
In view of Remark \ref{Rem3-new} and the assumption that $|\phi| \leq |\omega|$, it is enough to show that
$|\omega|\,|1-\omega|<1-|\omega|$ which, by squaring and then simplifying,  is seen to be equivalent to the
inequality
$$|\omega|^4-2{\rm Re\,}\omega |\omega|^2 +2|\omega|-1<0.
$$
In order to verify the last inequality, we observe that
\beqq
|\omega|^4-2{\rm Re\,}\omega |\omega|^2 +2|\omega|-1 &\leq & |\omega|^4+2|\omega|^3+2|\omega|-1\\
& =&(|\omega|^2+1)(|\omega|^2+2|\omega|-1)
\eeqq
which is negative whenever $|\omega|^2+2|\omega|-1<0$, i.e., $|\omega|<\sqrt{2}-1$.
The desired result follows.
\epf
\beg
If $|\phi| \leq |\omega|<s_0$ or $|\omega|\leq |\phi|<s_0$ on $\D$, then $C_{\psi,\phi}$ preserves
$\mathcal{P}_{\alpha}$ for every $\alpha$ with $-1<\alpha<0$, where $s_0 ~(\approx 0.2648)$ is the unique positive root of the
polynomial $P(x)=2x^4+8x^3+12x^2-1$.
\eeg
\bpf
Without loss of generality, we assume that $|\phi| \leq |\omega|$. In view of Remark \ref{p and q} and the assumption that
$\alpha \in (-1,0)$, the inequality \eqref{mp4 eq1} can be rewritten as
$$
 2|(1- \alpha)\,\overline{\omega}\,(\omega-1)+2i\alpha\, {\rm Im\,}\omega|\,|\phi|+(1- \alpha)(|\omega|^2
 +|1-\omega|^2|\phi|^2-1)-\alpha(1-|\omega|^2)(1-|\phi|^2)<0.
$$
By setting $\|\omega\|=A$ and $\|\phi\|=B$ (so that $B\leq A$), it sufficies to check that
$$
 2(1-\alpha)A(A+1)B-4\alpha AB+(1- \alpha)(A^2 +(1+A)^2B^2-1)- \alpha <0,
$$
which is equivalent to
$$
-\alpha[(A+B+AB)^2+4AB]+(A+B+AB)^2-1<0.
$$
Since $B\leq A$ and  $\alpha \in (-1,0)$, the last inequality holds if
$$
(2A+A^2)^2+4A^2]+(2A+A^2)^2-1 =2A^4+8A^3+12A^2-1<0.
$$
Clearly, the function $P(x)=2x^4+8x^3+12x^2-1$  is
strictly increasing on $(0,\infty)$ and thus,  $P(x)<0$ for $0\leq x< s_0$, where $s_0$
is the unique positive root of $P(x)$. The desired result follows.
\epf

\section{Fixed points}\label{MP4Sec6}
In this section, we discuss the fixed points of weighted composition operators. It is time to recall
a result due to Yu-Qing Chen \cite[Theorem 2.1]{Fixed points}. The modern way of writing it
is as follows:
\bprop
Let $X$ be a metrizable topological vector space and $C$ be a convex compact subset of $X$.
Then, every continuous mapping $T:C\rightarrow C$ has a fixed point in $C$.
\eprop

We set $X=\mathcal{H}(\mathbb{D})$, $C=\mathcal{P}_{\alpha}$, $T=C_{\psi,\phi}$ and observe that
every weighted composition operator on $\mathcal{P}_{\alpha}$ has a fixed point.
Indeed, one has something more than this as we can see below.

\bthm
Let $\phi$, $\psi$, $\omega$ be as before and $\phi$ is not a rotation. Suppose that $C_{\psi, \phi}$
is a self-map of $\mathcal{P}_{\alpha}$.  Then, $C_{\psi, \phi}$ has a
unique fixed point which can be obtained by iterating $C_{\psi, \phi}$ for any
$f \in \mathcal{P}_{\alpha}$. Further more, if $\phi$ is an inner function, then the fixed point
is the constant function $1$.
\ethm

\bthm
Let $\phi$, $\psi$, $\omega$ be as before and $\phi$ is  a rotation. Suppose that $C_{\psi, \phi}$
is a self-map of $\mathcal{P}_{\alpha}$ and $F$ denotes the set of all fixed points of $C_{\psi, \phi}$.
Then, there are three distinct cases:
\bee
\item If $\phi(z)\equiv z$, then $F=\mathcal{P}_{\alpha}$.
\item If $\phi(z)\equiv \lambda z$ and $\lambda^n\neq 1$ for every $n\in \N$, then $F=\{1\}.$
\item If $\phi(z)\equiv \lambda z$ and $\lambda^n=1$ for some $n>1$, then
$$ F=\{ f: f(z)=g(z^n)\mbox{~ for some ~} g\in\mathcal{P}_{\alpha} \}.
$$

\eee
\ethm

The proofs of these two theorems follow from the lines of the proofs of the corresponding
results of Section $4$ of \cite{arxiv W-CO}. Moreover, the key tools for the proofs are
from Section $6.1$ of Shapiro's book \cite{Shapiro:Book}. So we omit the details.

\subsection*{Acknowledgement}
The first author thanks the Council of Scientific and Industrial Research (CSIR), India,
for providing financial support in the form of a SPM Fellowship to carry out this research.
The second  author is currently at ISI Chennai Centre, Chennai, India.

\end{document}